\documentclass{gtart}
\usepackage{amssymb}
\let\Bbb\mathbb

\input gtmonout
\volumenumber{2}
\volumeyear{1999}
\volumename{Proceedings of the Kirbyfest}
\pagenumbers{135}{156}
\papernumber{8}
\received{12 January 1999}\revised{10 June 1999}
\published{17 November 1999}

\newcommand{\SO}{\mbox{\rm SO}}

\newcommand{\SL}{\mbox{\rm SL}}
\newcommand{\lcm}{\mbox{\rm lcm}}

\newcommand{\Spin}{\mbox{\scriptsize\rm Spin}} 
\newcommand{\Om}{\Omega_k^{\Spin}(B\pi )} 
\newcommand{\om}{\Omega_5^{\Spin}}
\newcommand{\omk}{\Omega_k^{\Spin}}
\newcommand{\Pos}{\mbox{\rm Pos}_k(\pi )} 
\newtheorem{thm}{Theorem}
\newtheorem{prop}[thm]{Proposition}

\newtheorem{cor}[thm]{Corollary}
\newtheorem{defn}[thm]{Definition}
\begin{document}
\title{Almost Linear Actions by Finite Groups on $S^{2n-1}$} 
\asciititle{Almost Linear Actions by Finite Groups on S^{2n-1}} 

\author{Hansj\"org Geiges\\Charles B Thomas}
\coverauthors{Hansj\noexpand\"org Geiges\\Charles B Thomas}
\asciiauthors{Hansjorg Geiges and Charles B Thomas}

\address{Mathematisch Instituut, Universiteit Leiden\\
Postbus 9512, 2300 RA Leiden, The Netherlands\\\smallskip
\\DPMMS, University of Cambridge\\
16 Mill Lane, Cambridge CB2 1SB, UK}
\asciiaddress{Mathematisch Instituut, Universiteit Leiden\\
Postbus 9512, 2300 RA Leiden, The Netherlands\\
\\DPMMS, University of Cambridge\\
16 Mill Lane, Cambridge CB2 1SB, UK}

\email{geiges@math.leidenuniv.nl\\C.B.Thomas@dpmms.cam.ac.uk}

\begin{abstract}
A free action of a finite group on an odd-dimensional sphere is said
to be {\sl almost linear} if the action restricted to each cyclic or
2--hyperelementary subgroup is conjugate to a free linear action. We
begin this survey paper by reviewing the status of almost linear
actions on the 3--sphere. We then discuss almost linear actions on
higher-dimensional spheres, paying special attention to the groups
$\SL_2(p)$, and relate such actions to surgery invariants. Finally, we
discuss geometric structures on space forms or, more generally, on
manifolds whose fundamental group has periodic cohomology. The
geometric structures considered here are contact structures and
Riemannian metrics with certain curvature properties.
\end{abstract}

\asciiabstract{A free action of a finite group on an odd-dimensional
sphere is said to be almost linear if the action restricted to
each cyclic or 2-hyperelementary subgroup is conjugate to a free
linear action. We begin this survey paper by reviewing the status of
almost linear actions on the 3-sphere. We then discuss almost linear
actions on higher-dimensional spheres, paying special attention to the
groups SL_2(p), and relate such actions to surgery
invariants. Finally, we discuss geometric structures on space forms
or, more generally, on manifolds whose fundamental group has periodic
cohomology. The geometric structures considered here are contact
structures and Riemannian metrics with certain curvature properties.}



\makeshorttitle
By the mid 1970s the existence problem for free actions by finite groups on $S^{2n-1}$, $n\geq 2$ (the topological spherical space form problem) had been solved. Classification has proved to be much harder, and there are residual problems in dimension three. Ib Madsen has proposed a definition of almost linear, applicable to a group $\pi$ which, although it does not itself act freely and linearly, may admit representations in $\mbox{\rm Top}^+(S^{2n-1})$ which are conjugate to representations in $\SO_{2n}$ when restricted to cyclic and 2--hyperelementary subgroups. This turns out to be the largest class of subgroups on which it is reasonable to impose a linearity condition. In his paper~\cite{mads78} Madsen conjectured that for all finite groups $\pi$ acting freely on $S^{2n-1}$ for some~$n$, there would exist almost linear actions, and gave a proof of this fact for certain $p$--hyperelementary groups of odd order. We provide further evidence for the truth of this conjecture by examining the groups $\SL_2(p)$. 

Given the development of surgery compatible with a geometric structure, it seems reasonable to give a survey of what is known about geometric spherical space forms, ie,\ space forms carrying an additional geometric structure (the structures considered here are positive scalar curvature metrics and contact structures, and we briefly touch on positive Ricci curvature metrics). In the geometric setting it is useful to weaken the condition of `almost linear' to `Sylow linear', that is to actions which admit linearization when restricted to a representative family of Sylow subgroups.

In dimension three the linearization of an arbitrary free $\pi$--action is part of the Thurston geometrization programme, but one which has proved resistant to his methods. Enough is known to show that the situation is very different from that in higher dimensions. Thus we have the omnibus theorem:

\medskip
{\bf Theorem}\qua(a)\qua {\sl Let $\pi$ be a finite solvable group acting freely and almost linearly on~$S^3$. Then the orbit space $S^3/\pi$ is Seifert fibred, and the action of $\pi$ is conjugate to a free linear action.}

(b)\qua {\sl Let $\pi$ be the binary icosahedral group $I^*\cong\SL_2(5)$. If $f\co\pi\rightarrow \mbox{\rm Top}^+(S^3)$ describes a free almost linear action, and $f_0$ describes the standard free linear action, then the join $f\ast f_0$ of the two actions, which is an action on $S^7$, is conjugate to a free linear action.} 
\medskip

We outline the proof of (a) in Section~1 below; part (b) is a special case of the construction given for $\SL_2(p)$ in Section~3. Note in passing that the only non-solvable finite groups which are candidates to act freely on $S^3$ are the direct products $C_u\times I^*$ with $\gcd (u,30)=1$. It is highly likely that the surgical argument for (b) extends to them.

Thus 3--dimensional topology is sufficiently rigid to remove the difference between linear and almost linear actions, at least when the group $\pi$ is solvable. In higher dimensions this fails, although there is some evidence that in dimension five almost linearity plus compatibility with a fixed geometric structure on $S^5$ may imply linearity, or at least that $\pi$ belongs to the class of groups known to have fixed point free representations in~$\SO_6$. This is a question which we intend to take up elsewhere, cf\ also our remarks in~\cite{geth98a}. 

In the present survey we cover the following material: \begin{enumerate}
\item Status of almost linear actions on~$S^3$. \item Almost linear actions on $S^{2n-1}$, $n\geq 2$, and surgery invariants. \item Special case of the finite algebraic groups $\SL_2(p)$. \item Geometric structures on space forms. \end{enumerate}
We would like the reader to become aware of the way in which almost linear actions draw on methods and results from a wide range of mathematics, including geometric methods in dimension three, advanced number theory, classical and not so classical differential geometry. In Section~5 we close by proposing some workpoints.

\section{Linearization in dimension three} %
We first list the (non-trivial) solvable finite groups which can act freely and linearly on the standard 3--dimensional sphere~$S^3$. 

\[ \begin{array}{rlcl}
(1) & D_{2^t(2v+1)}' & = & \left\{ a,b:\: a^{2v+1}=b^{2^t}=1,\, b^{-1}ab =a^{-1}\right\} ,t\geq 2,\; v\geq 1 \\ [2mm] (2) & Q_{4n} & = & \left\{ a,b:\; a^{2n}=1,\, a^n=b^2,\, b^{-1}ab=a^{-1} \right\} ,n\geq 2 \\ [2mm]
(3a) & \multicolumn{3}{l}{\mbox{\rm The generalized binary tetrahedral groups of order $8\cdot 3^v$}} \\ [2mm]
& T_v^* & = & \left\{ x,a,b:\; x^{3^v}=1,\, \langle a,b\rangle = Q_8,\, x^{-1}ax=b,\right. \\ [2mm]
& & & \left. \;\;\; x^{-1} bx=ab\right\} , v\geq 1\\ [2mm] (3b) & \multicolumn{3}{l}{\mbox{\rm The binary octahedral group of order 48}} \\ [2mm]
& O_1^* & = & \left\{ r,a,x,b:\; \langle r,a,b\rangle \cong Q_{16} ,\, \langle x,a,b\rangle = T_1^*,\, r^{-1}xr=x^{-1}\right\}\\ [2mm] (4) & \multicolumn{3}{l}{\mbox{\rm Direct products of any of the preceding groups (or the trivial group)}}\\
& \multicolumn{3}{l}{\mbox{\rm with a cyclic group $C_u$ of coprime order.}} \end{array} \]

For the details, which really amount to an exercise in representation theory, see \cite{miln57}, \cite{wolf67}, or~\cite{cbt86}. 

When referring to specific elements of these groups later on, we always adhere to the notation above.
The generator of $C_u$ will be denoted by~$y$. When speaking of a group in class 1, 2, or 3, we mean to include direct products with $C_u$. Observe that the groups in classes 1 and 2 have a cyclic subgroup of index two, generated by $ab^2y$ and $ay$, respectively.

In his 1957 paper, Milnor emphasized the importance of another family of groups with presentation
\[ \begin{array}{ccl}
Q(8n,k,l) & = & \left\{ z,a,b:\; \langle a,b\rangle =Q_{8n},\, z^{kl}=1,\, 
a^{-1}za=z^r,\right. \\
& & \left. \;\;\; b^{-1}zb=z^{-1}\right\} \;\mbox{\rm with} 
\; r\equiv -1\,\mbox{\rm mod}
\, k,\, r\equiv 1 \,\mbox{\rm mod} \, l. \end{array} \]
Again using representation theory~\cite[Section 8.2]{serr77}, one can see that there are free linear actions on spheres of dimension $8m-1$, and subject to delicate arithmetic conditions on $n,k,l$ (for example $n=7$, $k=809$, $l=1$, see \cite[pages~238 et seq.]{mads83}) the obstructions to the existence
of free topological actions in dimension $8m+3$ ($m\geq 1$) all vanish. Whether these exotic actions restrict to linear actions on a suitable class of subgroups is a very interesting question which deserves further investigation. We refer the reader to \cite{mads83} for the details of what is known.

In dimension three the existence of free actions by $Q(8n,k,l)$ depends on the existence of free non-linearizable cyclic group actions. This is explained in the survey article \cite{rubi95}, see also~\cite{cbt88}. A two-step covering space argument shows that if (at worst) each covering space whose fundamental group has order the product of two distinct odd primes is homeomorphic to a lens space, then no free $Q$--action can exist. The same pattern of argument also shows that in dimension three there is little difference between linear actions of a group $\pi$ and actions which are linear on each element of~$\pi$. The known geometric results are as follows:

\subsection*{$C_2$--actions}
\begin{thm}[Myers~\cite{myer81}]
{\sl A free involution on a $3$--dimensional lens space is conjugate to a free linear action.}
\end{thm}

Both the method and the result generalize G.R.~Livesay's theorem for free involutions on~$S^3$ \cite{live60}, in which the argument depends on finding an invariant {\em unknotted} circle~$S^1$. The action can then be linearized both in a tubular neighbourhood $S^1\times D^2$ and in its complement.
\subsection*{$C_3$--actions}
\begin{thm}[Rubinstein~\cite{rubi79a}]
{\sl If the cyclic group $C_3$ of order $3$ acts freely on ${\Bbb R}P^3$, then the orbit space is homeomorphic to a lens space.} \end{thm}

Let ${\cal M}$ be the class of orientable, closed, irreducible 3--manifolds with finite fundamental group, containing an embedded Klein bottle~$K^2$. This class constitutes the collection of Seifert fibred manifolds with $S^2$ as orbit surface and at most three exceptional fibres of multiplicity $2,2,n$ ($n\geq 1$), excluding $S^2\times S^1$, see~\cite{rubi79}. In particular, all manifolds in $\cal M$ are linear quotients of~$S^3$. 

\begin{thm}[Rubinstein~\cite{rubi79a}]
{\sl Let $C_3$ act freely on $M\in {\cal M}$. Then either $M/C_3$ is in $\cal M$, or $M/C_3$ is Seifert fibred over $S^2$ with three exceptional fibres of multiplicities $2,3,3$.}
\end{thm}

\begin{thm}[Rubinstein~\cite{rubi79a}]
{\sl Let the dihedral group $D_{2\cdot 3^s}$ act freely on a manifold $M\in {\cal M}$ in such a way that some element of order $2$ has an orbit space in $\cal M$. Then either the orbit space of the dihedral action belongs to $\cal M$ or is Seifert fibred over $S^2$ with three exceptional fibres of multiplicities $2,3,4$.}
\end{thm}

Using induction and the analogue of Theorem~3 for $C_2$--actions (see~\cite{rubi79}), one can extend Theorem~4 to dihedral groups of order $2^t3^s$ with $t>1$.

In the proof of all three theorems the technique is to consider an embedded ${\Bbb R}P^2$ or $K^2$ and its translates under $C_3$ or $D_{2\cdot 3^s}$. Rubinstein then simplifies the pattern of intersections far enough to show that the orbit space is as described. Several years ago Kirby pointed out to one of the authors that Rubinstein's original proof can be simplified to some extent by the use of minimal surfaces.

For the remainder of this section we discuss actions which satisfy the following restriction.

\medskip
{\bf Assumption}\qua ($C_{\mbox{\scriptsize\rm odd}}$)\qua Let $\pi$ be a finite group acting freely on $S^3$. If the element $a\in\pi$ has odd order, assume that the action restricted to the cyclic subgroup $\langle a\rangle$ generated by $a$ is conjugate to a free linear action.

\begin{thm}
{\sl If $\pi$ satisfies assumption $(C_{\mbox{\scriptsize\rm odd}})$, and $\pi$ is
one of the groups of type 1--4 listed above, then the action of $\pi$ is conjugate in $\mbox{\rm Top}^+(S^3)$ to a free linear action.} \end{thm}

{\bf Proof}\qua (i)\qua If $\pi$ is cyclic, repeated application of Myers' Theorem~1 gives the result.

(ii)\qua If $\pi$ is isomorphic to one of the groups in classes 1 or 2 (in the list at the beginning of this section), then $\pi$ contains a cyclic subgroup of index two, and Theorem~1 plus (i) apply. 

(iii)\qua Now consider the case that $\pi$ belongs to class 3a, ie, \[ \pi\cong C_u\times T_v^*\;\;\mbox{\rm with}\;\;\gcd (u,6)=1.\] The normal subgroup of index three generated by $\{ y,x^3,a,b\}$ is isomorphic to $C_{3^{v-1}u}\times Q_8$, and by step~(ii) the corresponding regular 3--fold covering space is the quotient of a free linear action, and hence Seifert fibred. By the classification of Seifert manifolds with a given fundamental group \cite[6.2]{orli72}, the quotient in question belongs to $\cal M$, hence Theorem~3 applies.

(iv)\qua Finally, if $\pi$ is modelled on the binary octahedral group $O_1^*$, it can be written as an extension
\[ C_u\times Q_8 \rightarrowtail \pi \twoheadrightarrow D_6.\] A 2--Sylow subgroup of $\pi$ is isomorphic to $Q_{16}$, and by~\cite{rubi79} the quotient of $S^3$ under the action of $Q_8$ and $Q_{16}$ is in~$\cal M$. Thus Theorem~4 applies.\endproof

Using the same technique (see~\cite{cbt88}) we can prove the following: 

\begin{cor}
{\sl There is no free action of $Q(8n,k,l)$ on $S^3$ satisfying assumption $(C_{\mbox{\scriptsize\rm odd}})$.} \end{cor}

{\bf Proof}\qua Write $Q(8p,q)$ for $Q(8p,q,1)$. Observe that $Q(8n,k,l)$ contains a subgroup $Q(8p,q)$, where either $p=2$ ($n$ even) or $p,q$ are distinct odd primes ($n$ odd), and we have extensions \[ C_q\rightarrowtail Q(16,q)\twoheadrightarrow Q_{16}\] and
\[ C_{pq}\rightarrowtail Q(8p,q)\twoheadrightarrow Q_8,\] respectively. The structural homomorphism $\varphi\co Q_{16}\rightarrow C_q$ (resp.\ $\varphi\co Q_8\rightarrow C_{pq}$) has kernel of order~8 (resp.\ of order~2).

Clearly it suffices to prove the corollary for $Q(8p,q)$. Arguing by contradiction, we assume that we are given a free action of $Q(8p,q)$ on $S^3$ satisfying $(C_{\mbox{\scriptsize\rm odd}})$. We have inclusion maps of subgroups of index~2,
\[ C_{4q}\hookrightarrow C_q\times Q_8\hookrightarrow Q(16,q)\;\;\; (p=2)\] and
\[ C_{2pq}\hookrightarrow Q_{4pq}\hookrightarrow Q(8p,q)\;\;\; 
(p\;\mbox{\rm odd}).\]
By the assumption $(C_{\mbox{\scriptsize\rm odd}})$ and Theorem~5, the fourfold covering space of $S^3/Q(8p,q)$ with fundamental group $C_{2pq}$ is homeomorphic to a lens space. A further application of Theorem~1 shows that the intermediate covering spaces are Seifert fibred, and as above we conclude from \cite{orli72} that they belong to $\cal M$. The simpler version of Theorem~3 for $C_2$--actions \cite[Theorem~8]{rubi79} shows that the orbit space under the second $C_2$--action is Seifert fibred and hence of linear type, which is a contradiction. \endproof 

%
%
\section{Free almost linear actions}
The groups $\pi$ considered in Section~1 are examples of {\sl periodic groups} (groups with periodic cohomology), that is, finite groups whose odd order Sylow subgroups are cyclic and whose 2--Sylow subgroups are either cyclic or generalized quaternionic, see~\cite[VI.9]{brow82}. In the case when the order of $\pi$ is odd, the group is necessarily {\sl metacyclic} (the extension of a cyclic group by another cyclic group) with presentation \[ \left\{ a,b:\; a^m=b^k=1,\, b^{-1}ab=a^r,\, \gcd ((r-1)k,m)=1,\, 
r^k\equiv 1\,\mbox{\rm mod}\, m\right\} ,\] see \cite[Theorem~5.4.1]{wolf67}. We recall that a group is said to be $p$--{\sl hyperelementary} if it is an extension of a cyclic normal subgroup by a $p$--group of coprime order. Subgroups of this kind, particularly when $p=2$, play a very important role in the existence and classification of topological spherical space forms. There is one significant difference between the cases $p=$ odd and $p=2$; in the former there may be no linear space forms, but in the latter representation theory always supplies examples. However, these may not be of the best possible dimension predicted by cohomology, see~\cite{wall78}. 

As an approximation of free linear actions, I~Madsen~\cite{mads78} proposed the following useful definition. Let ${\cal H}(\pi )$ denote the set of (conjugacy classes of) cyclic and 2--hyperelementary subgroups $\gamma$ of~$\pi$.

\begin{defn}
{\rm A free action of $\pi$ on $S^{2n-1}$ is said to be {\sl almost linear} if the action restricted to each $\gamma\in {\cal H}(\pi )$ is conjugate in $\mbox{\rm Top}^+(S^{2n-1})$ to a free linear action.} \end{defn}

We note that the set ${\cal H}(\pi )$ contains the Sylow subgroups $\pi_p$, so cohomological properties of the space form $M(\pi )=S^{2n-1}/\pi$, such as its normal invariant and smoothability, are detected at the level of the covering spaces $M(\gamma )$. In view of the geometric applications in Section~4 below it also makes sense to introduce the weaker notion of Sylow linearity.

\begin{defn}
{\rm A free action of $\pi$ on $S^{2n-1}$ is said to be {\sl Sylow linear} if the action restricted to each Sylow subgroup $\pi_p$ (for $p$ dividing the order $|\pi |$ of $\pi$) is conjugate in $\mbox{\rm Top}^+(S^{2n-1})$ to a free linear action.} \end{defn}

Even when the order of $\pi$ is odd (in which case all Sylow subgroups are cyclic), this is a strictly weaker notion, since almost linearity requires linearity for cyclic subgroups of composite order. There are examples of fake lens spaces with fundamental group $C_{pq}$ satisfying the Sylow linear condition, see \cite[Example~12.14]{miln66}.
Our aim in this section is to outline a proof (due to Madsen~\cite{mads78}) of the fact that in dimension greater than or equal to five, the strong condition of almost linearity does not imply linearity. We restrict attention to groups of odd order. The main steps of Madsen's construction are as described presently.

One key fact not available at the time of \cite{mads78} is the triviality of the group $SK_1({\Bbb Z}\pi )=\ker (K_1({\Bbb Z}\pi )\rightarrow K_1({\Bbb Q}\pi ))$ for periodic groups $\pi$ of odd order \cite[Theorem~14.2 (i)]{oliv88}. For finite groups, $SK_1( {\Bbb Z}\pi )$ equals the torsion subgroup of the Whitehead group $Wh(\pi )$, and if $SK_1({\Bbb Z}\pi )=0$, then Whitehead torsion is detected by Reidemeister torsion, weakly simple homotopy type is the same as simple homotopy type, and the surgery groups $L_*^s$ coincide with the intermediate surgery groups $L_*'$. For a non-technical account of surgery theory cf~\cite[Part~I]{wein94}.

\medskip
{\bf A}\qua Find a space type $X(\pi )$ (given by a finitely dominated $CW$ complex with an isomorphism $\pi_1(X)\rightarrow \pi$, a homotopy equivalence $\widetilde{X}\rightarrow S^{2n-1}$, and vanishing finiteness obstruction) whose covers $X(\gamma )$ have the simple homotopy type of given lens spaces, where $\gamma$ ranges over the Sylow subgroups of $\pi$ or, more optimistically, over all elements of ${\cal H}(\pi )$. The simple homotopy type of $X(\pi )$ is determined by the Reidemeister torsion $\Delta\in Wh({\Bbb Q}\pi )$, and the formal properties of $Wh$ guarantee, in the specific application below, that this simple homotopy type is uniquely determined by that of the covers.

\medskip
{\bf B}\qua Write $s(X(\pi ))$ for the topological structure set of $X(\pi )$, ie,\ equivalence classes of manifolds with a simple homotopy equivalence to $X(\pi )$. The surgery exact sequence reads \[ *\rightarrow s(X(\pi ))/L_{2n}^s(\pi )\rightarrow [X(\pi ),G/TOP] \rightarrow L^s_{2n-1}(\pi )=0, \]
where $L^s_{2n-1}(\pi )=0$ holds for any group of odd order~\cite{wall76}. Thus surgery is possible on all topological normal invariants in $[X(\pi ),G/TOP]$, and in this way we find a topological space form $M(\pi )$. We want to choose the normal invariant of $M(\pi )$ in such a way that it lifts to the normal invariants of given lens spaces.

\medskip
{\bf C}\qua Examine the action of $L_{2n}^s(\pi )$ on $s(M)$. To be more specific, it is necessary to know that $M(\pi )$ can be moved on its $L_{2n}^s(\pi )$--orbit in such a way that the covering spaces become homeomorphic to the given lens spaces.

\medskip
Steps {\bf A}, {\bf B}, and {\bf C} can certainly be carried out for a metacyclic group of odd order if we assume that $k=p^t$, ie,\ that $\pi$ is $p$--hyperelementary, and that $m$ is a power of a distinct prime~$q$. In order to rule out the uninteresting case when $\pi$ itself can act freely and linearly, we also assume that the action of $\langle b\rangle =C_{p^t}$ on $\langle a\rangle =C_m$ is
faithful.

Let $\chi_m$ and $\chi_p$ be free one-dimensional complex characters of the cyclic groups $C_m$ and $C_{p^t}$, respectively, and consider the representation spaces
\[ V_m=\bigoplus_{i=0}^{p^t-1} \chi_m^{r^i}\] and
\[ V_p=p \bigoplus_{i=0}^{p^{t-1}-1} \chi_p^{pi+1},\] both of real dimension $2n=2p^t$. The following theorem describes free non-linear but Sylow linear actions (hence almost linear actions for $m,k$ prime).

\begin{thm}[Madsen \cite{mads78}]
\label{thm:madsen}
{\sl Let $\pi$ be a metacyclic $p$--hyper\-elementary group with $p$ odd, and the order of the complementary normal subgroup $C_m$ also of odd prime power order. If the action of $C_{p^t}$ on $C_m$ is faithful, there exists a unique $($non-linear$)$ spherical space form $M^{2n-1}(\pi )$ with covering spaces $M(C_{p^t})$ and $M(C_m)$ homeomorphic to the lens spaces defined by $V_p$ and $V_m$, respectively.} \end{thm}

{\bf Remarks}\qua (1)\qua Non-linearity of $M(\pi )$ is a consequence of the fact that $\pi$ does not satisfy the $pq$--conditions, ie,\ not all subgroups of order $pq$ are cyclic. 

(2)\qua Madsen proved uniqueness up to weak $s$--cobordism; thanks to $SK_1({\Bbb Z}\pi )=0$
this can now be strengthened to uniqueness up to homeomorphism. 

(3)\qua Since the space form $M(\pi )$ obtained in the theorem is Sylow linear, the obstructions to smoothing $M(\pi )$ vanish when lifted to a Sylow covering space. As indicated before, this implies that $M(\pi )$ itself is smoothable (by the formal properties of the functor $[-,TOP/O]$), and the corresponding uniqueness statement still holds. Indeed, Madsen worked directly in the smooth category. 
\medskip

Here are some of the ideas in the proof of this theorem. As to {\bf B}, by Sullivan's theorem (cf~\cite{mami79},\cite[Section 2]{lama79}), we have \[ [M(\pi ),G/TOP]=\widetilde{KO}(M(\pi )),\] since $|\pi |$ is odd, and by \cite[Corollary~1.9]{mads78} \[ \widetilde{KO}(M(\pi ))\cong\widetilde{KO}(M(C_k))\oplus 
\widetilde{KO}(M(C_m))^{C_k}.\]
This decomposition already shows, up to the problem of $C_k$--invariance, that the normal invariant of $M(\pi )$ can be chosen to lift to the normal invariants for given lens spaces for $C_m$ and $C_k$, and that this choice is unique.

What concerns {\bf C}, we have a pull-back diagram \[ \begin{array}{ccc}
L^s_{2n}(\pi ) & \longrightarrow & L^s_{2n}(C_k) \\ \mbox{\Large $\downarrow$} & & \mbox{\Large $\downarrow$} \\ L_{2n}^s(C_m)^{C_k} & \longrightarrow & L_{2n}^s(1), \end{array} \]
where the maps to $L_{2n}^s(1)$ are surjective. We need to find $x\in L_{2n}^s(\pi )$ restricting to preassigned elements $x_m\in L_{2n}^s (C_m)$ and $x_k\in L_{2n}^s(C_k)$. One can arrange that $x_m$ and $x_k$ map to zero in $L_{2n}^s(1)$, so the existence of the desired $x$ is obvious from the diagram, provided one can show that $x_m$ is $C_k$--invariant. To achieve this, one reinterprets $x_m$ as an element in the complex character ring $R(C_m)$ via a homomorphism $\mbox{\rm sign}\co L_{2n}'(C_m) \rightarrow R(C_m)$ defined using the (multi)\-signature of the $\pi$--action. 
From the choices already made, one can deduce that $\mbox{\rm sign} (x_m)$ is 
$C_k$--invariant, and for $|\pi |$ odd the homomorphism sign is injective, which implies the invariance of~$x_m$.

It is interesting to contrast the argument just given with the earlier and geometrically more intuitive construction of T~Petrie \cite{petr71} using Brieskorn varieties. Consider the Brieskorn $(2p-1)$--manifold $\Sigma_{m,p}=\Sigma_{m,p}(\epsilon ,\eta )\subset {\Bbb C}^{p+1}$ defined by the equations
\[ z_0^m+...+z_{p-1}^m+z_p^l=\epsilon ,\] \[ \sum_{j=0}^p |z_j|^2=\eta .\]
Here $\epsilon$ and $\eta$ are positive real numbers, $m$ is an odd integer, and $l=p^s$ for some positive exponent~$s$. Let $\pi$ be the metacyclic group we have been considering, but with $k$ taken to be an odd prime number $p$ and $r$ to be a primitive $p$th root of~1 modulo~$m$. 

We have the short exact sequence
\[ 1\rightarrow C_m \stackrel{i}{\longrightarrow} \pi \stackrel{j}{\longrightarrow} C_p\rightarrow 1,\] where $C_m=\langle a\rangle$ and $C_p=\langle b\rangle$. Let $V_m$ and $V_p$ be free complex one-dimensional representations of $C_m$ and $C_p$, respectively (defined by primitive roots of unity $\xi_m,\xi_p$). Then $V=i_*V_m\oplus j^*V_p$ is a complex $(p+1)$--dimensional representation of~$\pi$. The action of $a$ multiplies the coordinates $z_i$ by $\xi_m^{r^i}$; the action of $b$ is a cyclic permutation on $(z_0,...,z_{p-1})$ and multiplies $z_p$ by~$\xi_p$. Then one can check directly that the action of $\pi$ on $\Sigma_{m,p}(\epsilon ,\eta )$ is free for suitable values of $\epsilon ,\eta$. This can be motivated by representation theory, where the intersection of the $\epsilon$--hypersurface with the $\eta$--sphere can be interpreted in terms of virtual representations. 

However, we are now faced with a surgical problem, since $H_{p-1}( \Sigma_{m,p};{\Bbb Z})$ is a finite group annihilated by $k=p^s$. The generators of this group can be removed by surgery (Petrie does this explicitly), and the surgery is effective because the relevant surgery obstruction group $L_{2p-1}^s(\pi )$ is trivial. In this way we obtain a free smooth action on some homotopy sphere in a dimension where the relevant group of homotopy spheres bounding parallelizable manifolds has order~1 or~2. A trick due to J~Morgan \cite[Theorem~5.3]{petr71} then allows us to take this homotopy sphere to be standard.

Each of the constructions we have outlined has its merits. But from the point of view of almost linear actions, the first is to be preferred, because its description of the $C_m$-- and $C_k$--coverings is so explicit. Since Petrie's construction requires surgery to remove the torsion middle homology group, it is not obvious that the intermediate coverings are lens spaces, because we do not know whether the signature of the resulting manifold is of linear type. The Reidemeister torsion is less of a problem, since it is calculated over the rational group ring ${\Bbb Q}\pi$, and hence does not `see' $H_{p-1}(\Sigma_{m,p}; {\Bbb Z})$.

\medskip
 {\bf Remark}\qua Note that the numerical restrictions on the two constructions are complementary. The first requires that the order of $\pi$ equal the product of two prime powers (though this can almost certainly be weakened); the second allows the order of $\langle a\rangle$ to be composite, but at the price of restricting that of $\langle b\rangle$ to be prime. %
\section{Special linear groups $\SL_2(p)$} %
In this section we consider the groups $\pi =\SL_2(p)$ of $2\times 2$ matrices with determinant 1 over the field ${\Bbb F}_p$, $p\geq5$ prime. Our basic references are \cite{lama79} and \cite{cbt97}. For a general value of $p$, the main facts from \cite{lama79} are the following:

\medskip
{\bf A}\qua Denote by $s(M(\pi ))$ the set of homeomorphism types of manifolds with a simple homotopy equivalence to some reference space form $M(\pi )$.
The latter exists by \cite{mtw76} and in the lowest possible dimension $4n-1=\lcm (4,p-1) -1$. The subgroup $L_0^s(1)\subset L_0^s(\pi )$ acts trivially on $s(M(\pi ))$, and the reduced group $\widetilde{L}_0^s (\pi )$ acts freely. The cardinality of the orbit space $s(M(\pi ))/\widetilde{L}_0^s(\pi )$ equals $p^s(p^2-1)^{n-1}$ with $s=4n/(p-1)-1$. This can be deduced from the surgery exact sequence \[ *\rightarrow s(M(\pi ))/\widetilde{L}_0^s(\pi )\rightarrow [ M(\pi ),G/TOP]\stackrel{\lambda}{\longrightarrow} L_3^s(\pi ),\] by observing that $\lambda$ vanishes identically (although $L_3^s(\pi )$ need not be trivial), and decomposing $[M(\pi ),G/TOP]$ into its local parts. The odd part, as in the previous section, comes from the real $K$--theory of $M(\pi )$, the 2--primary part from ordinary cohomology.

\medskip
{\bf B}\qua The simple homotopy type of $M(\pi )$ is determined by the Reidemeister torsion $\Delta (M)$, for once again we can appeal to a result of Oliver \cite[Theorem~14.15]{oliv88} that $SK_1({\Bbb Z}\pi )=0$ for the groups under consideration. Furthermore we have two classical surgery invariants: 

 (1)\qua The $\rho$--invariant
\[ \rho \co s(M(\pi ))\longrightarrow \widetilde{RO}(\pi )\otimes {\Bbb Q},\] where $RO(\pi )$ is the real representation ring of $\pi$ and \[ \widetilde{RO}(\pi )=\mbox{\rm coker}(i_*\co RO(1)\rightarrow RO(\pi )),\] that is, $\widetilde{RO}(\pi )$ is the quotient of $RO(\pi )$ under the ideal generated by the regular representation. 

 (2)\qua The 2--local normal invariant \[ \nu_{(2)}\co s(M(\pi ))\longrightarrow \bigoplus_i H^{4i}(M; 
{\Bbb Z}_{(2)}).\]
One then has the result that the number of topological manifolds contained in a fixed simple homotopy type $\Delta$ and with given invariants $\rho$ and $\nu_{(2)}$ is bounded by the cardinality of $\mbox{\rm Tor}\, L_0^s(\pi )$. This follows from the surgery exact sequence above, together with the fact that $\rho$ determines the odd-primary part of the normal invariant in $[M(\pi ),G/TOP]$ and detects the action of $\widetilde{L}_0^s(\pi )/\mbox{\rm Tor}$ (because of the formula relating it to the multisignature \[ \mbox{\rm sign}\co \widetilde{L}_0^s(\pi )\longrightarrow 
\widetilde{RO}(\pi )\otimes {\Bbb Q},\]
given by
\[ \rho (\alpha\cdot \{ f\})=\rho (\{ f\} )-\mbox{\rm sign} (\alpha ), \;\; \alpha\in \widetilde{L}_0^s(\pi ),\; \{ f\} \in s(M(\pi )),\] and the injectivity of sign on the torsion-free part of~$\widetilde{L}_0^s(\pi )$).

We now want to discuss almost linear actions of $\pi =\SL_2(p)$. To do this, we begin with a more detailed description of its subgroup structure, cf\ \cite{lama79}. There are three basic cyclic subgroups $C_{(1)},C_{(2)},C_{(3)}$, whose conjugates exhaust~$\pi$: \[ \begin{array}{rcl}
C_{(1)} & = & \left\{ \pm\left( \begin{array}{cc}1 & \alpha \\ 
0 & 1
\end{array} \right) ,\alpha\in{\Bbb F}_p\right\}
\cong C_{2p}, \\ [5mm]
C_{(2)} & = & \left\{ \left( \begin{array}{cc} \beta & 0\\ 
0 & \beta^{-1}
\end{array} \right) ,\beta\in {\Bbb F}_p^{\times}
\right\} \cong C_{p-1}, \\ [5mm]
C_{(3)} & = & \ker\left\{ \mbox{\rm norm homomorphism}\; 
N:{\Bbb F}_{p^2}^{\times}\rightarrow {\Bbb F}_p^{\times}, \gamma\mapsto \gamma\cdot \gamma^p \right\} \cong C_{p+1}.
\end{array} \]
In the last definition, ${\Bbb F}_{p^2}$ is regarded as a 2--dimensional vector space over ${\Bbb F}_p$, and then $\gamma\in {\Bbb F}_{p^2}^{\times}$ gives an element of $\mbox{\rm GL}_2(p)$ with determinant as described. 

Maximal 2--hyperelementary subgroups $\tau_1,\tau_2,\tau_3$ are given by \[ \begin{array}{ccccccccc}
1 & \rightarrow & C_p & \rightarrow & \tau_1 & \rightarrow & 
C_{2^l} & \rightarrow 1, \\
1 & \rightarrow & C_m & \rightarrow & \tau_2 & \rightarrow & 
C_4 & \rightarrow 1, \\
1 & \rightarrow & C_n & \rightarrow & \tau_3 & \rightarrow & 
Q_{2^{k+1}} & \rightarrow 1,
\end{array} \]
where $C_{2^l}$ is the 2--Sylow subgroup of ${\Bbb F}_p^{\times}$, $m$ and $n$ are odd, and $\{ 2m,2^kn\} = \{ p-1, p+1\}$. Then $\tau_2\subset N(C_{(3)})$ and $\tau_3\subset N(C_{(2)})$ if $p\equiv 1$ mod~4, otherwise $\tau_i\subset N(C_{(i)})$, $i=1,2$. Any other 2--hyperelementary subgroup
is contained in a conjugate of one of these models. Finally, by hyperelementary induction one obtains (see~\cite[Lemma~3.1]{lama79}) \[ L_*^s(\pi )=L_*^s(Q_{2^{k+1}})^{\mbox{\scriptsize\rm stable}} \oplus L_*^s(\tau_1)_{\mbox{\scriptsize\rm odd}}^{{\Bbb F}_p^{\times}} \oplus L_*^s(\tau_2)_{\mbox{\scriptsize\rm odd}} \oplus L_*^s(\tau_3)_{\mbox{\scriptsize\rm odd}}. \] Here
\[ L_*^s(\tau_i)=L_*^s(\tau_{i(2)})\oplus L_*^s (\tau_i)_{\mbox{\scriptsize\rm odd}} ,\] because each extension splits, and ${\Bbb F}_p^{\times}$ acts on $L_*^s(\tau_1)$ through its quotient by its 2--Sylow subgroup. 

Combining all this information, one can prove the following, see~\cite{cbt97}. For the purpose of this survey, we state the result in a slightly weaker form than there.

\begin{thm}
{\sl For suitable values of $n$ divisible by $4$ and $p-1$ there exists a free topological action of $\pi=\SL_2(p)$ on $S^{2n-1}$ such that torsion and normal invariants of this action are those of a linear action when restricted to the subgroups $\tau_2$ and $\tau_3$.} \end{thm}

{\bf Idea of the proof}\qua The existence of some free action is implied by the main theorem in~\cite{mtw76}. To be more precise, start with a fixed point free representation of $\pi_2\cong Q_{2^{k+1}}$ and choose compatible representations for $\tau_2$ and $\tau_3$. Next choose a Reidemeister torsion $\Delta$ which restricts to the given linear space form torsions for these two subgroups. To do this in the minimal dimension $p-2$ (for $p\equiv 1$ mod~4) resp.\ $2p-3$ (for $p\equiv 3$ mod~4) requires a careful choice of these representations; a cheaper version is available if we allow the dimension to be larger, see~\cite[Remark~2.12]{mads78}.

The odd part of the normal invariant lies in real $K$--theory by Sullivan's theorem quoted in Section~2, and the decomposition of the $\widetilde{KO}$--groups obtained by hyperelementary induction is analogous to that of the $L$--groups given above. Thus the odd part of the normal invariant of $M(\pi )$ can be chosen to restrict to the given linear ones for $\tau_2$ and $\tau_3$. We can do the same for the 2--local normal invariant by an application of~\cite[Lemma~3.3]{mtw76}.\endproof

To obtain an action that actually restricts to a linear action, we need to complete the step analogous to {\bf C} in Section~2. Hence the problem is that the preassigned elements $\alpha_2\in L_0^s(\tau_2)$ and $\alpha_3 \in L_0^s(\tau_3)$ given by the choice of $M(\pi )$ and the linear representations for $\tau_2$ and $\tau_3$ have to come from an element $\alpha\in L_0^s(\pi )$, so the decomposition formula for $L_0^s(\pi )$ given above shows that the 2--primary parts of $\alpha_2$ and $\alpha_3$ have to be stable (ie,\ invariant under the action of $N(Q_{2^{k+1}})$). The $\rho$--invariants of $M(\tau_2)$ and $M(\tau_3)$ are stable by construction, so if the representations of $\tau_2$ and $\tau_3$ are chosen in such a way that the corresponding $\rho$--invariants are stable, then the remarks in {\bf B}(2) show that the torsion-free part of $\alpha_2$ and $\alpha_3$ is stable. Finally, from the formula for $L_0^s(\pi )$ we have that $\mbox{\rm Tor}L_o^s(\pi )$ restricts onto $\mbox{\rm Tor}L_0^s(\tau_2 ) \oplus \mbox{\rm Tor}L_0^s(\tau_3 )$, since $L_0^s$ is torsion-free for cyclic and quaternion 2--groups by~\cite{wall76}. 

\vspace {2mm}

A heuristic explanation of this theorem is given by the fact that lifting the inclusion representation $\SL_2(p)\rightarrow \SL_2( \overline{{\Bbb F}}_p)$ to characteristic zero gives a virtual representation of virtual dimension 2 over ${\Bbb C}$ which restricts to a positive fixed point free representation of each of the cyclic groups $C_{(2)}$ and $C_{(3)}$.

Sharper results can be obtained by restricting the prime~$p$, see \cite{cbt97} for the case $p=2^t+1$; similar considerations are possible for $p=2^t-1$. The case when $p=5$ is most interesting, since it has implications in dimension three. The following theorem gives part (b) of the theorem in the introduction.

\begin{thm}
{\sl If $\SL_2(5)$ acts freely on $S^3$ in such a way that the action restricted to a cyclic subgroup of order $5$ is conjugate to a free linear action, then the $\SL_2(5)$--action embeds in a free linear action on~$S^7$.}
\end{thm}

{\bf Sketch proof}\qua The discussion above shows that the maximal 2--hyperelemen\-tary subgroups are $\tau_1\cong Q_{20}$, $\tau_2\cong Q_{12}$ and $\tau_3\cong Q_8$. If $f\co\SL_2(5)\rightarrow \mbox{\rm Top}^+(S^3)$ is the representation defining the action, then $f|\tau_i$ ($i=1,2,3$) is conjugate to a free linear action. For $\tau_2$ and $\tau_3$ this follows immediately from Theorems~1 and~2; for $\tau_1$ we use the assumption plus two applications of Theorem~1. Furthermore, each linear orbit space is unique up to homeomorphism. Since the Reidemeister torsion $\Delta$ and signature $\rho$ are detected on hyperelementary subgroups, it follows that these two algebraic invariants for $f$ coincide with those of a free linear action~$f_0$. (Representation theory shows that $\SL_2(5)$ has two fixed point free irreducible representations of dimension~2, switched by the unique outer automorphism of the group.) The mod~2 part of the normal invariant is similarly fixed since it is 2--local and detected on~$\tau_3$.

Now let $\SL_2(5)$ act on $S^7=S^3\ast S^3$, the join of two copies of~$S^3$, via $f$ on the first factor and $f_0$ on the second. All three invariants $\Delta$, $\rho$ and $\nu_{(2)}$ are compatible with this construction. This is clear from the definitions for $\Delta$ and $\nu_{(2)}$, for $\rho$ see~\cite[page~185]{wall70}. It follows that the join of the two actions is linear on $S^7$; note that because the primes involved are small there is no 2--torsion in~$L_0^s(\SL_2(5))$. \endproof

In anticipation of surgery in the presence of a geometric structure, let us consider the smoothing of the topological space forms obtained. Since the map $\mbox{\rm BO}\rightarrow\mbox{\rm BTop}$ is a map of infinite loop spaces, $k(-)=[-,\mbox{\rm B(Top/O)}]$ is a respectable cohomology theory with transfer, from which it follows that \[ k(M(\pi ))\longrightarrow\bigoplus_{l|\, |\pi |} k(M(\pi_l))\] is injective. Therefore any obstruction to smoothing is detected at the level of Sylow subgroups, ie,\ Sylow linear space forms are smoothable. In the case at hand this implies that $M(\SL_2(p))$ will certainly be smoothable if the covers $M(\tau_i)$, $i=1,2,3$, are linear. 

In sufficiently low dimensions, it is probable that one can prove smoothability without specifying the structure of $M(\tau_1)$. Provided $M(\tau_2)$ and $M(\tau_3)$ are linear, the remaining obstruction to smoothability is detected in the cohomology of a $p$--Sylow subgroup (cyclic of order~$p$) with coefficients in $(\mbox{\rm Top/O})^* (\mbox{\rm point})$. By \cite[Theorem~5.18]{mami79} we have a splitting of $H$--spaces $\mbox{\rm G/O}_{(p)}\simeq \mbox{\rm BSO}_{(p)}\times \mbox{\rm coker}J_{(p)}$. This suggests that, having reduced from G to Top, further reduction from Top to O depends on $(\mbox{\rm coker}J)_* \subset\pi_*^s$. But the first $p$--torsion occurs in dimension $p(2p-2)-2$ \cite[Theorem~1.1.14]{rave86}, which is much higher than the expected minimum dimension $p-2$ or $2p-3$ of $M(\pi )$. This again fits in with the earlier observation about lifting from characteristic~$p$. %
\section{Geometric structures and periodic fundamental\nl groups} %

Starting with the work of Gromov--Lawson and Schoen--Yau on metrics of
positive scalar curvature, the methods of surgery, used to solve
topological problems as in Sections 2 and~3 above, have been shown to
be applicable in various geometric settings.

Recall that the Levi--Civita connection associated with a Riemannian
metric $g$ on a smooth manifold $M$ gives rise to a fourth order
covariant curvature tensor $R_{ijkl}$, satisfying various $(\pm 1)$
symmetric and Jacobi relations. A first contraction gives the
symmetric Ricci curvature tensor
\[ (\mbox{\rm Ric})_{jk}=g^{li}R_{ijkl},\] 
and a further contraction gives the scalar curvature 
\[ s=g^{jk}(\mbox{\rm Ric})_{jk}.\]
This scalar invariant can be thought of as measuring the average
sectional curvature at a point of~$M$; the Ricci curvature, viewed as
a quadratic form on tangent directions, measures the average sectional
curvature over 2--planes containing the given direction. In either
case, the geometrically interesting problem is that of the existence
of {\em positive} curvature metrics, there being no obstructions to
the existence of negative Ricci curvature metrics in dimension greater
than two~\cite{lohk94}.

Recall also that a {\sl contact structure} on an odd-dimensional manifold $M^{2n+1}$ is a maximally non-integrable hyperplane field $\xi\subset TM$, that is, defined locally as the kernel of a 1--form $\alpha$ satisfying $\alpha\wedge (d\alpha )^n\neq 0$. If $\xi$ is coorientable, it can be written as $\xi =\ker\alpha$ with a globally defined $\alpha$, which is then called a {\sl contact form}. 

These then are the three geometric structures in which we are interested. What concerns compatibility with surgery along an embedded sphere, the situation at the time of writing is the following: 

\begin{itemize}
\item Positive scalar curvature (psc) metrics: Surgery is possible along spheres of codimension~$\geq 3$~\cite{grla80}. \item Contact structures: Surgery is possible along isotropic spheres (ie,\ spheres tangent to the contact structure) of less than half the dimension, together with a framing restriction. Approximation by isotropic spheres is governed by an $h$--principle and can thus be controlled by topological conditions (such as the vanishing of certain characteristic classes), see \cite{elia90}, \cite{wein91}, \cite{geig97}, \cite{geth98}. \item Positive Ricci curvature (pRc) metrics: This is perhaps the most interesting, and certainly the least well understood of the three geometric structures under consideration. Here surgery is possible in dimension~$\geq 2$ and codimension~$\geq 3$, provided the metric in the neighbourhood of the surgered sphere is in a suitable standard form~\cite{wrai98}. The codimension bound alone is almost certainly too weak in general, and there is a putative codimension~$\geq 7$ condition, though there is no clear indication how Wraith's construction in~\cite{wrai98}, which is a refinement of one by Sha--Yang and Schoen--Yau, should gain in flexibility by allowing for higher codimension. For a discussion of these issues see~\cite{stol96}.
\end{itemize}

We begin by considering psc metrics. Following \cite{kwsc90}, let \[ \mbox{\rm Pos}_k(\pi )\subset\Om\]
be those classes represented by
$(V\rightarrow B\pi ,\sigma )$, where $V$ is a psc manifold with spin structure~$\sigma$. Then $\mbox{\rm Pos}_k(\pi )$ is a subgroup, and by Gromov--Lawson surgery one can easily prove the following statement: 

\begin{prop}
{\sl Let $f\co M\rightarrow B\pi $ represent a class in $\mbox{\rm Pos}_k(\pi )$, $k\geq 5$. If $f$ is $2$--connected, then $M$ admits a psc metric.} \end{prop}

The subgroups $\Pos$ have good naturality properties: 

\begin{itemize}
\item covariant naturality: If $\varphi\co \pi\rightarrow \pi'$ is a group homomorphism, composition with $\varphi$ defines an induced map \[ (B\varphi )_*\co\Om\longrightarrow \omk (B\pi ').\] 
From the definitions, $(B\varphi )_*$ sends $\Pos$ to $\mbox{\rm Pos}_k(\pi ')$. 
\item contravariant naturality: If $\varphi \co\pi\hookrightarrow \pi '$ is an inclusion, there is a transfer homomorphism \[ (B\varphi )^!\co\omk (B\pi ')\longrightarrow \Om \] defined as follows: Starting from a given map $f'\co V'\rightarrow B\pi '$, we obtain a pullback diagram
\[ \begin{array}{ccc}
V & \stackrel{f}{\longrightarrow} & B\pi \\ \mbox{\Large $\downarrow$} & & \;\;\mbox{\Large $\downarrow$} B\varphi \\ V' & \stackrel{f'}{\longrightarrow} & B\pi '. \end{array} \]
The given spin structure $\sigma '$ lifts to a spin structure $\sigma$ on $V$, and since a psc metric also lifts to a covering space, $(B\varphi )^!$ maps $\mbox{\rm Pos}_k(\pi ')$ to~$\Pos$. \end{itemize}

\begin{prop}[Kwasik--Schultz \cite{kwsc90}] {\sl With $j_p\co \pi_p\rightarrow\pi$ denoting inclusion of a $p$--Sylow subgroup of $\pi$, we have that $\alpha\in\Pos$ if and only if $(Bj_p)^! \alpha\in\mbox{\rm Pos}_k(\pi_p)$ for all primes $p$ dividing the order of~$\pi$.}
\end{prop}

Modulo an elementary lemma on
Noetherian modules this follows from the fact that after localization at~$p$, the composition $(Bj_p)_*(Bj_p)^!$ is an isomorphism; see \cite{kwsc90} for the details.

This proposition allows to reduce the exitence problem for psc metrics on manifolds of dimension $\geq 5$ and fundamental group $\pi$ to the corresponding problem for the $\pi_p$. For periodic fundamental groups this has been used to show that such a manifold admits a psc metric if and only if its universal cover does (Gromov--Lawson--Rosenberg conjecture). Kwasik--Schultz~\cite{kwsc90} dealt with groups of odd order (where all $\pi_p$ are cyclic); by the corresponding study of generalized quaternion groups, the general case was settled by Botvinnik--Gilkey--Stolz~\cite{bgs97}.

In dimension $2n+1=5$ the codimension $\geq 3$ and dimension $\leq n$ conditions coincide. In the case of contact manifolds, the relevant $h$--principle tells us that attaching 1--spheres can always be taken to be isotropic, and the same holds for 2--spheres on which the first Chern class of the contact structure vanishes. The framing condition is either stably satisfied (for $S^1$) or empty (for $S^2$). Thus, by mimicing the argument of Kwasik--Schultz, we obtain the following theorem, see~\cite{geth}. 

\begin{thm}
\label{thm:geth}
{\sl Let $\pi$ be a finite periodic group of odd order not divisible
by 9. Then every closed $5$--dimensional spin manifold $M$ with fundamental group isomorphic to $\pi$ admits a contact structure.} \end{thm}

Both this theorem and the results of Kwasik--Schultz and Botvinnik et al.\ apply to the space forms constructed in Section~2. So far as psc metrics are concerned, the same holds for the $\SL_2(p)$--examples in dimension $p-2$ or $2p-3$. As yet we can say nothing about the existence of contact structures, since for $p\geq 7$ the dimension is too high. 

By combining the methods used to prove Theorem~\ref{thm:geth} with Madsen's construction, one can give explicit generators of the spin bordism groups $\Om$ (with $\pi$ as in Theorem~\ref{thm:madsen}) in terms of space forms that carry contact structures (with analogous statements in the psc case). For instance, as explained in~\cite{geth}, for the metacyclic group $\pi =D_{q^r,3}$ with normal subgroup $C_{q^r}$ and quotient group $C_3$, one can show by spectral sequence arguments that $\om (BD_{q^r,3})$ is a cyclic group of order $9p^r$, generated by a 5--dimensional spherical space form carrying a contact structure. 

What concerns pRc metrics on non-linear spherical space forms, the best hope to obtain these might be to use Petrie's construction described in Section~2. By the work of Hern\'andez--Andrade~\cite{hern75}, suitably chosen $\Sigma_{m,p}$ admit $\pi$--invariant pRc metrics, so the problem reduces to being able to perform pRc surgery on $\Sigma_{m,p}/\pi$. Unfortunately, while the description of $\Sigma_{m,p}/\pi$ is explicit, that of the surgery necessary to obtain a space form is not, so a stronger surgery theorem than Wraith's is required. Compare \cite[pages~279/80]{kwsc90} and, for related issues in the psc case, \cite{kwsc96}.

The strongest curvature condition one can impose is that the sectional curvature be positive. To the best of our knowledge, there are no examples of non-linear free actions of a finite group on an odd-dimensional sphere with an invariant metric of strictly positive sectional curvature. In a recent paper, however, Grove and Ziller~\cite{grzi} show the existence of infinitely many non-isometric metrics of non-negative sectional curvature on each of the four homotopy projective 5--spaces (ie,\ quotients of $S^5$ under some free ${\Bbb Z}_2$--action; cf~\cite{geth98} for a thorough discussion of these spaces in the context of contact geometry). They start from a description of $S^5$ as a Brieskorn variety carrying an action by $\SO_2\times \SO_3$ of cohomogeneity one (ie,\ an action with one-dimensional orbit space). This action descends to a cohomogeneity one action on $S^5/{\Bbb Z}_2$, with ${\Bbb Z}_2$ generated by $-\mbox{\rm id}\in\SO_2$, and all four homotopy projective 5--spaces arise in this way. Moreover, the singular orbits of this action are of codimension two, and Grove--Ziller give a general gluing construction for producing metrics of non-negative sectional curvature on precisely such cohomogeneity one manifolds with codimension two singular orbits.
\section{Some workpoints}

(1)\qua Continue the work of I~Madsen and J~Milgram on the groups $Q(8n,k,l)$. As is pointed out in \cite{mads83} the basic groups to consider are the $Q(8p,q)$, the semi-direct products of $C_{pq}$ and $Q_8$, where $p$ and $q$ are distinct odd primes, and the structural map from $Q_8$ into $\mbox{\rm Aut}(C_{pq})$ has kernel of order~2. In a covering space argument reminiscent of that sketched in Section~1, Madsen detects the surgery obstruction for existence in the quaternionic subgroups of~$Q$, but does not consider classification. Various questions suggest themselves, for example:
\begin{itemize}
\item[(a)] What values are possible for the Reidemeister torsion $\Delta$ and signature $\rho$?
\item[(b)] Do there exist almost linear (or Sylow linear) actions by $Q(8p,q)$ on $S^{8m+3}$ with $m>0$?
\end{itemize}

(2)\qua Complete the construction of free almost linear actions for periodic groups of odd order. First extend Theorem~\ref{thm:madsen} to arbitrary $p$--hyperelementary groups, and then to metacyclic groups $\pi$ for which the second factor $k$ of the order is composite. 

(3)\qua Calculate the signatures of the actions constructed by Petrie. 

(4)\qua Give a more geometric construction of the free actions of $\SL_2(p)$ constructed in Section~3. Since the identity representation in characteristic~$p$ lifts to the difference of a $(p+1)$-- and $(p-1)$--dimensional representation in characteristic zero, there is the possibility of using the intersection of two hypersurfaces picked in the manner of Petrie. This hardly helps with the prime $p$ with its period $p-1$, and so the next step may be to construct a free action on the product $S^3\times S^{2p-1}$ rather than on a single sphere. 

(5)\qua In Section~3 we used surgery to embed an almost linear action of $I^*=\SL_2(5)$ on $S^3$ into the standard linear action on~$S^7$. Is it possible to desuspend this action in such a way as to show that the original action on $S^3$ was also standard? It may be worthwile remembering that $S^7$ carries an $I^*$--invariant pRc metric, and that the action is linear on any invariant subsphere $S^3$ on which the induced metric also satisfies the pRc condition by a result of Hamilton~\cite{hami82}, which uses heat flow methods to construct constant positive curvature metrics on 3--manifolds with pRc metrics. 

(6)\qua Improve the conditions for contact surgery to obtain results on the existence
of contact structures parallel to those of Kwasik--Schultz in dimensions~$\geq 7$. This will involve an analysis of framing conditions, but the major obstacle will be to have sufficient control over the positions of the attaching spheres above the middle dimension. 

(7)\qua Study contact structures on fake lens spaces. This problem can be posed independently of the previous one, since we can start with free $C_p$--actions on Brieskorn varieties studied by Hirzebruch and others. So long as further modi\-fications take place along embedded spheres of less than half the dimension, existing methods will suffice. %

\Addresses\recd

\end{document}